\documentclass[12pt]{article}
\usepackage{amssymb}
\usepackage{latexsym}
\usepackage[all]{xy}
\newcounter{mytab}
\newenvironment{mytab}[1]{
\refstepcounter{mytab}
\begin{center}
Table \arabic{mytab}. #1
\end{center}
}{}
\newtheorem{THM}{Theorem}[section]
\newtheorem{PROP}{Proposition}[section]

\newcommand{\LegD}[2]{
( \,{#1}\,| \, {#2}\, )_{\mathbb{D}}
}
\newcommand{\Leg}[2]{
(\, {#1}\,| \, {#2}\, )
}
\begin{document}
%
%%%%%%%%%%%%%%%%%%%%%%%%%%%%%%%%%%%%%%%%%%%%%%%%%%%%%%%%%%%%%%%%
%  

\title{Modular reduction in abstract polytopes}
\author{B. Monson\thanks{Supported by NSERC of Canada Grant \# 4818} 
\hskip.08in and Egon Schulte\thanks{Supported by NSA-grant 
H98230-07-1-0005}}
\maketitle
\begin{abstract}
\noindent
The paper studies modular reduction techniques for abstract regular and chiral  
polytopes, with two purposes in mind:\  first, to survey the literature about 
modular reduction in polytopes; and second, to apply modular reduction, with 
moduli given by primes in $\mathbb{Z}[\tau]$ (with $\tau$ the golden ratio), 
to construct new regular $4$-polytopes of hyperbolic types 
$\{3,5,3\}$ and $\{5,3,5\}$ with automorphism groups given by finite 
orthogonal groups.



\end{abstract}

\section{Introduction}
\label{intro}

Polytopes and their symmetry have inspired mathematicians since antiquity. 
In the past three decades, a modern abstract theory of polytopes has emerged 
featuring an attractive interplay of mathematical areas, including geometry, 
combinatorics, group theory, graph theory and topology (see \cite{arp}). 
Abstract polytopes share many properties with ordinary convex polytopes but a 
priori are not embedded in the geometry of an ambient space.

The present paper studies modular reduction techniques for regular and chiral 
abstract polytopes. Modular reduction has proved to be a powerful tool in the 
construction and analysis of new classes of polytopes. Our paper serves two 
purposes:\ first, to survey  the literature about modular reduction techniques 
in polytopes; and second, to apply modular reduction, with moduli given by 
primes in the ring of integers $\mathbb{Z}[\tau]$ of the quadratic number 
field $\mathbb{Q}(\sqrt{5})$, to construct new infinite classes of regular 
$4$-polytopes of hyperbolic types $\{3,5,3\}$ and $\{5,3,5\}$ with 
automorphism groups given by orthogonal groups over finite fields.

Last, but not least, as a birthday greeting, we wish to acknowledge the many 
contributions by our friend and colleague Ted Bisztriczky in convexity, 
polytope theory and combinatorial geometry (see also \cite{bis}).
 
\section{Basic notions and methods}
\label{absgps} 

For general background material on abstract polytopes we refer the reader 
to \cite[Chs. 2,3]{arp}. Here we just review some basic terminology.

An (\emph{abstract\/}) \emph{polytope of rank\/} $n$, or an 
\emph{$n$-polytope\/}, is a partially ordered set $\mathcal{P}$ with a 
strictly monotone rank function having range $\{-1,0, \ldots, n\}$. An 
element of rank $j$ is called a \emph{$j$-face\/} of $\mathcal{P}$, and a 
face of rank $0$, $1$ or $n-1$ is a \emph{vertex\/}, \emph{edge\/} or 
\emph{facet\/}, respectively. The maximal chains, or \emph{flags}, of 
$\mathcal{P}$ all contain exactly $n + 2$ faces, including a unique least 
face $F_{-1}$ (of rank $-1$) and a unique greatest face $F_n$ (of rank $n$). 
Two flags are said to be \emph{adjacent} ($i$-\emph{adjacent}) if they differ 
in a single face (just their $i$-face, respectively). Then $\mathcal{P}$ is 
required to be \emph{strongly flag-connected} (see \cite[Ch.2]{arp}). Finally, $\mathcal{P}$ has the following homogeneity 
property:\  whenever $F \leq G$, with $F$ a $(j-1)$-face and $G$ a 
$(j+1)$-face for some $j$, then there are exactly two $j$-faces 
H with $F \leq H \leq G$.
  
Whenever $F \leq G$ are faces of ranks $j \leq k$ in $\mathcal{P}$, the 
\textit{section}
$G/F := \{ H \in \mathcal{P}\, | \, F \leq H \leq G \}$ is thus 
a ($k-j-1$)-polytope in its own right. In particular, we can
identify $F$ with $F/F_{-1}$. Moreover, we call 
$F_{n}/F$ the \emph{co-face at\/} $F$, or the 
\emph{vertex-figure at\/} $F$ if $F$ is a vertex.

Our interest is primarily 
in regular or chiral polytopes. A polytope $\mathcal{P}$ is \emph{regular\/} 
if its \textit{automorphism group} $\Gamma(\mathcal{P})$ is transitive on the flags of $\mathcal{P}$, and $\mathcal{P}$ is \emph{chiral\/} if $\Gamma(\mathcal{P})$ has two flag 
orbits such that adjacent flags are always in distinct orbits. 

For a regular polytope $\mathcal{P}$, the group $\Gamma(\mathcal{P})$ 
is generated by $n$ involutions $\rho_0,\ldots,\rho_{n-1}$, where 
$\rho_i$ maps a fixed, or \emph{base\/}, flag $\Phi$ to the flag 
$\Phi^i$, $i$-adjacent to $\Phi$. These generators satisfy (at least) 
the standard Coxeter-type relations 
\begin{equation}
\label{standrel}
(\rho_i \rho_j)^{p_{ij}} = \epsilon \;\; \textrm{ for } i,j=0, \ldots,n-1, 
\end{equation}
where $p_{ii}=1$, $p_{ji} = p_{ij} =: p_{i+1}$ if $j=i+1$, and 
$p_{ij}=2$ otherwise; thus the underlying Coxeter diagram is a string 
diagram. Note that
$p_{j} := p_{j-1,j} \in \{2,3,\ldots,\infty\}$ for $j = 1,\ldots,n-1$. 
These numbers
 determine 
the (\emph{Schl\"afli}) \emph{type} $\{p_{1},\ldots,p_{n-1}\}$ of 
$\mathcal{P}$. Moreover, the following {\em intersection condition\/} holds:
\begin{equation}
\label{intprop}
\langle \rho_i \mid i \in I \rangle \cap \langle \rho_i \mid i \in J \rangle 
= \langle \rho_i \mid i \in {I \cap J} \rangle 
\qquad (I,J \subseteq \{0,1,\ldots,n-1\}). 
\end{equation}
The rotations 
$\sigma_{i}:=\rho_{i}\rho_{i-1}\;(i=1,\ldots,n-1)$
generate the {\em rotation subgroup\/} $\Gamma^{+}(\mathcal{P})$ of 
$\Gamma(\mathcal{P})$, which is of index at most~$2$. We call 
$\mathcal{P}$ {\em directly
regular\/} if this index is $2$.

A group $\Gamma = \langle \rho_0, \ldots, \rho_{n-1} \rangle$ whose 
generators satisfy (\ref{standrel}) and (\ref{intprop}), is called a  
\emph{string C-group\/}; here, the ``C'' stands for ``Coxeter'', though 
not every C-group is a Coxeter group. These string C-groups are precisely 
the automorphism groups of regular polytopes, 
since, in a natural way, such a polytope can be uniquely reconstructed from
$\Gamma$ (see \cite[\S 2E]{arp}). Therefore, we   often identify a regular 
polytope with its automorphism (string C-) group. 

The group $\Gamma(\mathcal{P})$ of a chiral polytope $\mathcal{P}$ is 
generated by $n-1$ elements $\sigma_1,\ldots,\sigma_{n-1}$, which again are 
associated with a base flag $\Phi=\{F_{1},F_0,\ldots,F_{n}\}$, such that 
$\sigma_i$ fixes all the faces in $\Phi \setminus \{F_{i-1},F_{i}\}$ and 
cyclically permutes (``rotates") consecutive $i$-faces of $P$ in the 
(polygonal) section $F_{i+1}/F_{i-2}$ of rank $2$. By replacing a 
generator by its inverse if need be, we can further require that, 
if $F_{i}'$ denotes the $i$-face of $\mathcal{P}$ with 
$F_{i-1}<F_{i}'<F_{i+1}$ and $F_{i}'\neq F_i$, then 
$\sigma_{i}(F_{i}')=F_i$. The resulting distinguished generators 
$\sigma_1,\ldots,\sigma_{n-1}$ of $\Gamma(\mathcal{P})$ then satisfy
relations
\begin{equation}
\label{chirrel}
\sigma_i^{p_i} = (\sigma_j\sigma_{j+1}\ldots\sigma_k)^{2} = \epsilon 
 \;\; \textrm{ for } i,j,k=1,\dots,n-1,  \textrm{ with } j<k, 
\end{equation}
where again the numbers $p_i$ determine the type 
$\{p_1,\ldots,p_{n-1}\}$ 
of $P$. Moreover, $\Gamma(\mathcal{P})$ and its generators satisfy a certain 
intersection condition resembling that for C-groups. Conversely, if $\Gamma$ 
is a group generated by $\sigma_1,\ldots,\sigma_{n-1}$ such that the 
relations (\ref{chirrel}) and the new intersection condition hold, then 
$\Gamma$ is the group of a chiral polytope, or the rotation subgroup 
for a directly regular polytope; here the 
polytope is regular if and only if $\Gamma$ admits an involutory 
automorphism $\rho$ such that 
$\rho(\sigma_1)=\sigma_{1}^{-1}$, $\rho(\sigma_2)=\sigma_{1}^{2}\sigma_2$, and 
$\rho(\sigma_j)=\sigma_j$ for $j\geq 2$. For a 
chiral polytope $\mathcal{P}$, the two flag orbits yield two sets of 
generators $\sigma_i$ which are not conjugate in $\Gamma(\mathcal{P})$; 
thus a chiral polytope 
occurs in two \emph{enantiomorphic} (mirror image) 
forms.

We now describe the basic idea of modular reduction. We begin with a linear 
group $G$ over a ring $\mathbb{D}$, choose an ideal $J$ of $\mathbb{D}$, and try 
to construct a polytope from the quotient group of $G$ obtained by viewing $G$ 
as a linear group over $\mathbb{D}/J$. For the latter step the main obstruction 
is typically the intersection condition for the resulting quotient group of $G$. 
In most applications, $G$ itself is already a string C-group 
(often a Coxeter group), $\mathbb{D}$ is the ring of integers in an 
algebraic number field, and $J$ is an ideal of $\mathbb{D}$.   

More precisely, let $\mathbb{D}$ be a commutative ring with identity $1$, 
let $V$ be a free module over $\mathbb{D}$ of rank $n$ with basis 
$b_{0},\ldots,b_{n-1}$, and let $G$ be a subgroup of the general 
linear group $GL_{n}(\mathbb{D})$ over $\mathbb{D}$, whose elements 
we may view as invertible linear transformations (module isomorphisms) of $V$. 
Now let $J$ be an ideal of $\mathbb{D}$. Then the natural ring epimorphism 
$\mathbb{D} \rightarrow \mathbb{D}/J$ defined by $a \rightarrow a+J$ induces a 
group epimorphism $G \rightarrow G^J$, where $G^J$ is a subgroup of 
$GL_{n}(\mathbb{D}/J)$, the \emph{modular reduction\/} of $G$; in other words, 
we obtain $G^J$ by simply viewing the matrix entries of the elements in $G$ as 
entries in $\mathbb{D}/J$. Then $G^J$ naturally acts on the free module $V^J$ 
over $\mathbb{D}/J$ of rank $n$ with basis $b_{0}^J,\ldots,b_{n-1}^J$. We often 
abuse notation by referring to the modular images of objects by the same name 
(such as $V$ or $b_{0},\ldots,b_{n-1}$, etc.); that is, we drop the 
superscript $J$.

Modular reduction is, of course, a natural idea and has been used in various
 ways to construct maps and polytopes \cite{monwei1,monwei2,wil1}. Here we
begin by applying the idea to crystallographic Coxeter groups.

\section{Polytopes from crystallographic Coxeter groups}
\label{polcry}

Let $\Gamma$ be an abstract {\em string\/} Coxeter group with generators 
$\rho_{0},\ldots,\rho_{n-1}$ and presentation as in (\ref{standrel}), where 
again $p_{ii} = 1$ and $p_{ij} = 2$ for $|i-j| \geq 2$. Let $V$ be real 
$n$-space, with basis $\alpha = \{a_{0},\ldots,a_{n-1}\}$ and symmetric 
bilinear form  $x \cdot y$ defined by
\begin{equation}
\label{stanfm}
a_i \cdot a_j := - 2 \cos{\textstyle \frac{\pi}{p_{ij}}},\;\;\;0
\leq i,j \leq n - 1.
\end{equation}
Let $R: \Gamma \rightarrow G$ be the (faithful) standard representation of 
$\Gamma$ in $V$, where $G = \langle r_{0},\ldots,r_{n-1} \rangle$ is the 
isometric reflection 
group on $V$ generated by the reflections with \textit{roots} $a_i$ 
(see \cite[\S 5.3--5.4]{humph}); thus, 
\[ r_i(x) = x - (x \cdot a_i)\, a_i \quad (i = 0,\ldots,n-1). \] 
Now let $m = 2m'$, where $m'$ is the lowest common multiple of all $p_{ij}$ 
which are finite.
Let $\xi$ be a primitive $m$-th root of unity, and let 
$\mathbb{D}:= \mathbb{Z} [ \xi ]$. Then, with respect to the basis $\alpha$ of $V$, 
the reflections $r_j$ are represented by matrices in $GL_n(\mathbb{D})$ so that 
we may view $G$ as a subgroup of $GL_n(\mathbb{D})$. 
(By \cite[Th. 21.13]{curein}, $\mathbb{D}$ is the ring of integers in the 
algebraic number field $\mathbb{Q} (\xi)$; and $\mathbb{D}$ has (finite) 
rank $\phi(m)$ as a $\mathbb{Z}$-module.) Now we can
{\em reduce $G$ mod $p$}, for any prime $p$, here allowing $p = 2$ 
\cite[ch. XII]{curein}. More precisely, suppose that $p$ is a rational 
prime and 
that $J$ is a maximal ideal in $\mathbb{D}$ with 
$p\mathbb{D} \subseteq J \subset \mathbb{D}$. Then 
$\mathbb{K}:= \mathbb{D}/J$ is a 
finite field of characteristic $p$, and reduction mod $p$ of $G$ is achieved 
by applying the natural epimorphism $\mathbb{D} \rightarrow \mathbb{K}$ to the 
matrix entries of the elements in $G$.  This then defines a representation 
$\kappa: G \rightarrow GL_n(\mathbb{K)}$ with image group 
$G^p := G^J= \kappa(G)$. (This construction is essentially independent of the 
choice of $J$.)  Note that $\kappa$ is faithful when $G$ is finite 
and $p \nmid |G|$ (in fact, $\kappa$ often is faithful even when $p$ divides 
$|G|$). In any case, $\ker\kappa$ is a $p$-subgroup of $G$. 

The modular reduction technique is considerably more straightforward for 
{\em crystallographic\/} Coxeter groups, meaning that $G$ (or $\Gamma$) leaves 
invariant some lattice in $V$. A (string) Coxeter group $G$ is known to be 
crystallographic if and only if $p_{ij}=2$, $3$, $4$, $6$ or $\infty$ for 
all $i \neq j$. In fact, 
even allowing non-string diagrams,
it is also true that $G$ is crystallographic if 
and only if there is a \textit{basic system} 
$\beta = \{b_{0},\ldots,b_{n-1}\}$, with $b_i:=t_i a_i$ for certain 
$t_i > 0$, such that 
$m_{ij} := - t_{i}^{-1} (a_{i} \cdot a_{j})\, t_{j} \in \mathbb{Z}$  for
 $0 \leq i,j \leq n-1$. 
Then, for the rescaled roots $b_i$, we have
\begin{equation}
\label{reflnII}
r_i(b_j) = b_j + m_{ij} b_i ,  
\end{equation}
so that the generators $r_i$ are represented by integral matrices 
with respect to the basis $\beta$, and the corresponding {\em root lattice} 
${\displaystyle \oplus_{j} \mathbb{Z} b_j}$ actually is $G$-invariant. 
Thus we may take $\mathbb{D}=\mathbb{Z}$, and $G$ can be reduced modulo 
any integer $d\geq 2$. We focus primarily on the case when $p$ is an 
odd prime, but address some questions regarding composite moduli as well. 
For a composite modulus $d\geq 2$, the reduced group $G^d$ generally 
does not ``split'' according to the prime factorization of $d$, and its 
structure depends more heavily on the diagram $\Delta(G)$ used in the 
reduction process (two diagrams which modulo a prime are equivalent, 
may not be equivalent modulo a composite modulus). 
We now describe these diagrams.

The various possible basic systems $\{t_{i}a_{i}\}$ for a crystallographic 
Coxeter group $G$ can be represented by a diagram $\Delta(G)$ 
(see \cite[p. 415]{extf} or \cite{sqf}):\  for 
$0 \leq i,j \leq n-1$, node $i$ is labelled $2t_i^2\:(= b_{i}^2)$; and 
distinct nodes $i \neq j$ are joined by 
$\lambda_{ij} := \min\{m_{ij}, m_{ji}\}$ 
{\em unlabelled} branches.  Note that 
$\lambda_{ij} = \lambda_{ji} = 0$, $1$ or $2$, so the 
underlying graph is essentially that of the underlying 
Coxeter diagram $\Delta_c(G)$ of $G$, except that a 
mark $p_{ij} = \infty$ is indicated by a {\em doubled\/} 
branch in the case that $m_{ij} = m_{ji} =2$.  In 
Table~\ref{binartabII} we display the possible 
subdiagrams corresponding to the dihedral subgroups 
$\langle r_i, r_j\rangle$. For simplicity we have replaced 
the node labels $2t_i^2, 2t_j^2$ by $s,t$ or $s, ks \;(k=1,2,3,4)$, 
as appropriate. Note here that 
$m_{ij}m_{ji} = 4 \cos^2 \frac{\pi}{p_{ij}}$, so that, if 
$m_{ij} \geq m_{ji} \geq 1$, we have 
${t_{j}^2 / t_{i}^2} = {m_{ij} / m_{ji}} = 1,2,3,4\mbox{ (or }1)$ for 
$p_{ij} = 3,4,6,\infty$, respectively. 

\begin{center}
\begin{tabular}{c|c|c} 
Nodes $i, j$ & Parameters & $\lambda_{ij}$ \\ 
\hline
$ \stackrel{s}{\bullet}\;\;\;\;\stackrel{t}{\bullet}$ & $p_{ij} = 2$ &
	0 	 \\
\hline
$ \stackrel{s}{\bullet}\!\!\!-\!\!\!-\!\!\!-\!\!\!\!\stackrel{k s}{\bullet}$ 
	& $p_{ij} = 3, 4, 6, \infty$ & 1\\  
	& ($k = 1, 2, 3, 4$) &  \\ 
\hline
$ \stackrel{s}{\bullet}\!\!=\!\!=\!\!=\!\!\stackrel{s}{\bullet}$
	&$p_{ij} = \infty$ & 2 \\
	& ($k = 1$) & 
\end{tabular}
\end{center}
\begin{mytab}{Basic Systems for the Crystallographic 
Dihedral Groups $\langle r_i, r_j\rangle$.}
\label{binartabII}
\end{mytab}

The Gram matrix $B = [b_{ij}] := [b_i \cdot b_j]$ is easily computed 
from $\Delta(G)$, since 
$b_{ii} = 2 t_i^2$ is simply the label attached to node $i$, and
$b_{ij} = -\lambda_{ij}\max\{b_{ii}, b_{jj}\}/2$ for $i \neq j$. 
Moreover,
$m_{ij} = \lambda_{ij} \max\{1, b_{jj}/b_{ii}\}$ for $i \neq j$.

For a connected Coxeter diagram $\Delta_{c}(G)$, the corresponding group 
$G$ has to scale 
only finitely many basic systems $\beta$. Any such system is represented 
by an essentially 
unique diagram $\Delta (G)$, in which node labels form a set of relatively 
prime integers.   
For two such basic systems $\beta, \beta'$, we can convert from $\Delta(G)$ 
to $\Delta'(G)$ by consecutively adjusting the labels and branches of various 
pairs of adjacent nodes by operations of the following kind:\  inverting the 
ratio of the node labels; doubling a single branch and balancing its labels, 
or converting a double branch to a single branch, with ratio $4$, if the 
corresponding branch in $\Delta_{c}(G)$ is marked $\infty$.  Following 
these adjustments on pairs of nodes, we may finally have to rescale the 
entire set of labels to obtain a set of relatively prime integers. 

Then we may reduce $G$ modulo an odd prime $p$ to obtain a subgroup $G^p$ 
of $GL_n(\mathbb{Z}_p)$ generated by the modular images of the $r_i$'s.  
We call the prime $p$  \textit{generic} for $G$ if $p\geq 5$, or $p=3$ but 
no branch of $\Delta_{c}(G)$ is marked $6$. In the generic case, no node 
label of $\Delta(G)$ vanishes mod $p$ and a change in the underlying basic 
system for $G$ has the effect of merely conjugating $G^p$ in 
$GL_n(\mathbb{Z}_p)$. On the other hand, in the non-generic case, the group 
$G^p$ may depend essentially on the actual diagram $\Delta(G)$ used in the 
reduction mod $p$. For any odd prime $p$, we always find that 
$G^p = \langle r_0, \ldots , r_{n-1} \rangle$ is a subgroup of the orthogonal 
group $O( \mathbb{Z}_p^n)$ of isometries for the (possibly singular) symmetric 
bilinear form $x \cdot y$,  the latter being defined on $\mathbb{Z}_p^n$ by 
means of (the modular image of) the Gram matrix $B$; in particular, $r_{i}$ 
is the orthogonal reflection with {\em root\/} $b_i$ if $b_i^2 \neq 0$.  

Now recall from \cite[Thm. 3.1]{monsch1} that an irreducible group $G^p$ 
of the above sort, 
generated by $n \geq 3$ reflections, must necessarily be one of the following: 
\begin{itemize}
\item an orthogonal group $O(n,p,\epsilon) = O(V)$ or 
$O_j(n,p,\epsilon) = O_{j}(V)$, excluding the cases 
$O_1(3,3,0)$, $O_2(3,5,0)$, $O_2(5,3,0)$ (supposing for these three 
that 
$\mbox{disc}(V) \sim 1$), and also excluding the case $O_j(4,3,-1)$; or 
\item the reduction mod $p$ of one of the finite linear Coxeter groups of 
type $A_n$ 
($p\nmid n+1$), $B_n$, $D_n$, $E_6$ ($p \neq 3$), $E_7$, $E_8$, $F_4$, $H_3$ 
or $H_4$.  
\end{itemize}
We shall say in these two cases that $G^p$ is of \textit{orthogonal} or 
\textit{spherical type}, respectively, although there is some overlap for 
small primes. (The notation $a \sim b$ means that $a = t^2 b$ for some $t$ 
in the field.) Our description in \cite[Thm. 3.1]{monsch1} rests on the 
classification of the finite irreducible reflection groups over any field, 
obtained in {Zalesski\u\i} \& {Sere\v zkin} \cite{zal1} 
(see also \cite{kant,wag1,wag2,zal2}). For a non-singular space $V$, recall 
that $\epsilon = 0$ if $n$ is odd, and $\epsilon = 1$ or $-1$ if $n$ is even 
and the Witt index of $V$ is $n/2$ or $(n/2)-1$.  When $n$ is even, 
$\mbox{disc}(V) \sim (-1)^{n/2}$ if $\epsilon = 1$, and 
$\mbox{disc}(V) \sim (-1)^{n/2}\gamma$, with $\gamma$ a non-square, if 
$\epsilon = -1$. Moreover, $O_{1}(V)$ and $O_{2}(V)$ are the subgroups of 
$O(V)$ generated by the two distinct conjugacy classes of reflections, 
each characterized by the quadratic character of the spinor norm of its 
reflections, namely $1$ or $\gamma$, respectively. 

The above analysis sets up the stage for the construction of finite regular 
polytopes from  crystallographic string Coxeter groups $G$. Clearly, the 
generators $r_i$ of the reduced group $G^p$ satisfy (at least) the Coxeter-type
relations inherited from $G$. In most cases it also is more or less straightforward to 
determine the overall structure of $G^p$ by appealing to the above description 
of possible reflection groups. However, the main challenge is to determine 
when $G^p$ has the intersection property (\ref{intprop}) for its standard 
subgroups.  Here the outcome is often quite unpredictable, with the result 
depending on the underlying group $G$, the diagram $\Delta(G)$ used in the 
reduction, and the prime $p$. The reduced group $G^p$ may turn out to be a 
C-group for all primes, or for all primes in certain congruence classes, or 
for only finitely many primes (for example, only for $p=3$), or for no prime 
at all. Thus there is no theorem that covers all cases simultaneously.

Here we describe a quite general theorem (see \cite[Thm. 2.3]{monsch3}) 
which, in the three basic cases considered, typically enables us to determine 
(often after additional considerations) whether or not the reduced group
 $G^p$ is a C-group. For any $k,l \in \{0, \ldots n-1\}$ we let 
 $G_k^p := \langle r_j \,|\, j \neq k\rangle$  and 
 $G_{k,l}^p := \langle r_j \,|\, j \neq k,l \rangle$. 
 Then $G^p$ is a string C-group if and only if $G_{0}^p$ and 
 $G_{n-1}^p$ are string C-groups and 
 $G_{0}^p \cap G_{n-1}^p = G_{0,n-1}^p$. 
 (This enables an inductive attack on the problem.)
 Similarly we let $V_k$ and $V_{k,l}$ denote the subspace of 
 $V = \mathbb{Z}_p^n$ spanned by  the vectors $b_j$ with $j\neq k$ or 
 $j\neq k,l$, respectively. For a \textit{singular} subspace $W$ of a 
 non-singular space 
$V$, we also let $\widehat{O}(W)$ denote the subgroup of $O(W)$ 
consisting of those isometries which act trivially on the radical of $W$.

\begin{THM}
\label{subspcrit}  
Let $G = \langle r_0, \ldots , r_{n-1} \rangle$  be a crystallographic 
linear Coxeter group with string diagram. Suppose  that $n\geq 3$, that 
the prime $p$ is generic for $G$ and that there is a square among the 
labels of the nodes $1,\ldots,n-2$ of the diagram $\Delta (G)$ 
(this can be achieved by readjusting the node labels). For various subspaces 
$W$ of  $V$ we identify 
$O(W)$, $\widehat{O}(W)$, etc. with suitable subgroups of the pointwise 
stabilizer of $W^\perp$ in $O(V)$.

\noindent{\rm (a)}
Let the subspaces  $V_{0}$, $V_{n-1}$ and  
$V_{0,n-1}$ be non-singular, and let $G_0^p$, $G_{n-1}^p$ be of orthogonal 
type.  
 
{\rm (i)}   Then $G_{0}^p \cap G_{n-1}^p$ acts trivially on 
$V_{0,n-1}^{\perp}$ and
$ O_{1}(V_{0,n-1}) \leq G_{0}^p \cap G_{n-1}^p \leq  O (V_{0,n-1})$. 

{\rm (ii)} If $G_{0}^{p} = O(V_{0})$ and $G_{n-1}^{p} = O(V_{n-1})$,  then 
$ G_{0}^p \cap G_{n-1}^p = O (V_{0,n-1})$. 

{\rm (iii)} If either $G_{0}^{p} = O_1(V_{0})$  or 
$G_{n-1}^{p} = O_1(V_{n-1})$,  then 
$ G_{0}^p \cap G_{n-1}^p = O_{1}(V_{0,n-1})$.   

\noindent{\rm (b)}
Let $V$, $V_{0}$, $V_{n-1}$ be non-singular, let $V_{0,n-1}$ be singular 
(so that $n \geq 4$), 
and let $G_0^p, G_{n-1}^p$ be of orthogonal type.  

{\rm (i)}  Then  $G_{0}^p \cap G_{n-1}^p$ acts trivially on 
$V_{0,n-1}^{\perp}$, and
$\widehat{O}_{1}(V_{0,n-1}) \leq G_{0}^p \cap G_{n-1}^p \leq  
\widehat{O} (V_{0,n-1})$.

{\rm (ii)} If $G_{0}^{p} = O(V_{0})$ and $G_{n-1}^{p} = O(V_{n-1})$, then  
$ \widehat{O}(V_{0,n-1}) = G_{0}^p \cap G_{n-1}^p $. 

{\rm (iii)}  If either $G_{0}^{p} = O_1(V_{0})$ or  
$G_{n-1}^{p} = O_1(V_{n-1})$, then 
$ \widehat{O}_1(V_{0,n-1}) = G_{0}^p \cap G_{n-1}^p $. 

\noindent{\rm (c)}  Suppose $V , V_{0,n-1}$ are non-singular 
while at least one of  the subspaces $V_{0}$, $V_{n-1}$ is  singular.
Also suppose that $G_{0,n-1}^p$ is  of orthogonal type, with
 $G^p = O_1(V)$  when  $G_{0,n-1}^p = O_1(V_{0,n-1})$. Then 
$G_{0}^p \cap G_{n-1}^p  = G_{0,n-1}^p$.
\end{THM}

We briefly discuss the modular polytopes associated with some interesting classes of crystallographic string Coxeter groups. Recall that $[p_{1},p_{2},\ldots,p_{n-1}]$ denotes the Coxeter~group 
with a string Coxeter diagram on $n$ nodes and branches labelled 
$p_{1},p_{2},\ldots,p_{n-1}$. They are the automorphism groups of the 
universal regular polytopes
$\{p_{1},p_{2},\ldots,p_{d-1}\}$ (see \cite[\S 3D]{arp}).
\vskip.06in

\noindent {\it 3.1.\ Groups in which $G$ has spherical or Euclidean Type}

The modular reduction $G^p$ of any spherical or euclidean (crystallographic)
 group $G$ is a string C-group for any prime $p\geq 3$, and $G^p\cong G$ if 
 $G$ is spherical.  There are four kinds of (connected) spherical string 
 diagrams (up to duality), namely $A_n$, $B_n$,  $F_4$ and
 $I_2(6)$  (dihedral of order $12$). The 
 corresponding modular polytopes are isomorphic to the $n$-simplex, $n$-cube, 
$24$-cell and hexagon, respectively, and admit 
``modular realizations" in the 
 finite space $\mathbb{Z}_p^n$. For the euclidean groups $[4,3^{n-3},4]$ 
 ($n\geq 3$), $[3,4,3,3]$ and $[\infty]$, we obtain the regular toroids 
 $\{4,3^{n-3},4\}_{(p,0^{n-2})}$ of rank $n$ and $\{3,4,3,3\}_{(p,0,0,0)}$ of 
 rank $5$, and the regular $p$-gon $\{p\}$.

For a spherical group $G$, we have $G^d \simeq G$ for any modulus $d\geq 3$, 
and sometimes even for $d=2$ (here depending on choice of diagram).
For a euclidean group $G$ and $d\geq 3$ 
(again, sometimes for $d=2$ as well), the reduced group $G^d$ is the group of
a regular toroid of rank $n$, but now its type vector (suffix) depends on 
$\Delta(G)$ (in particular, on the parity of $n$) and, as well, 
on  the parity of $d$; the 
details are quite involved (see \cite{monschmod}).

\noindent {\it 3.2.\ Groups of ranks $3$ or $4$}

The groups of ranks $1$ or $2$ are subsumed by our discussion in the 
previous paragraph. For rank $3$, any group $G^p$ is a string C-group. For 
example, the hyperbolic group $G=[3,\infty]$ with diagram
\[ \stackrel{1}{\bullet}\frac{}{\;\;\;\;\;\;\;}
\stackrel{1}{\bullet}\frac{}{\;\;\;\;\;\;\;}
\stackrel{4}{\bullet}\;\;\;(\mbox{disc}(V) = -1), \]
yields a regular map of type $\{3,p\}$. For $p\geq 5$ we find that 
$G^p = O_1(3,p,0) \simeq PSL_2(\mathbb{Z}_p) \rtimes C_2$ of order 
$p (p^2 -1)$. When $p=3$, $5$ or $7$, respectively, we obtain the 
regular tetrahedron $\{3,3\}$, icosahedron $\{3,5\}$ and the Klein 
polyhedron $\{3,7\}_8$. This construction redescribes the family 
of regular maps discussed (and generalized) in \cite{mcmpl} or 
\cite{regmaps}.
 
The situation changes drastically for groups $G$ of higher ranks, 
with obstructions already occurring in rank $4$. We shall not attempt 
to fully summarize the findings of \cite{monsch2}, which settle all 
groups $G=[k,l,m]$ of rank $4$ (employing results such as 
Theorem~\ref{subspcrit}). Suffice it here to mention some possible 
scenarios illustrating that the outcome is often unpredictable. 

One possible scenario is that $G^p$ is a string C-group for all primes 
$p\geq 3$; this happens, for instance, if the subgroup $[k,l]$ or 
$[l,m]$ is spherical, if $l=\infty$, or if $[k,l]$ or $[l,m]$ is 
euclidean and $l=4$ or $6$. For example, from $[3,3,\infty]$ we 
obtain a regular $4$-polytope of type $\{3,3,p\}$ whose vertex-figures 
are isomorphic to the maps of type $\{3,p\}$ derived from $[3,\infty]$; 
its group $G^p$ is isomorphic to $S_5$ if $p=3$, $O_{1}(4,3,1)$ if 
$p\equiv 1 \bmod{4}$, and $O_{1}(4,3,-1)$ if $p\equiv 3 \bmod{4}$, 
$p\neq 3$. When $p=3$, $5$ or $7$, respectively, this gives the 
$4$-simplex $\{3,3,3\}$, the $600$-cell $\{3,3,5\}$, and the universal 
$4$-polytope $\{\{3,3\} , \{3,7\}_{8}\}$ first described in \cite{monw1}. 

On the other extreme, there are groups $G$ for which $G^p$ is a C-group for 
only finitely many primes $p$. For example, for $G=[6,3,6]$ only $G^3$ is a 
C-group. More important examples arise when the four subspaces $V$, $V_0$, 
$V_3$, $V_{0,3}$ are all non-singular. For example, when 
$G=[\infty,3,\infty]$ the reduced group $G^p$ is a C-group only for 
$p=3$, $5$ or $7$; from $G^3$ and $G^5$ we obtain the $4$-simplex 
$\{3,3,3\}$ and the regular star-polytope $\{5,3,\frac{5}{2}\}$, and from 
$G^7$ a (non-universal) $4$-polytope with facets and vertex-figures given 
by a dual pair of Klein maps $\{7,3\}_8$ and $\{3,7\}_8$.

Another possible scenario is that $G^p$ is a C-group only for certain 
congruence classes of primes. If $G=[6,3,\infty]$, then $G^p$ is a C-group 
if and only if $p=3$ or $p\equiv \pm 5 \bmod{12}$.

\noindent {\it 3.3.\ Groups of higher ranks}

The large number of crystallographic Coxeter groups $G$ of ranks 
$n\geq 5$ makes it difficult to 
fully enumerate the regular polytopes obtained by our method. 

For the group $[4,3,4,3]$ of rank $5$ we obtain regular $5$-polytopes with 
toroids $\{4,3,4\}_{(p,0,0)}$ as facets and $24$-cells $\{3,4,3\}$ as 
vertex-figures; their group $G^p$ is given by $O_1(5,p,0)$ if 
$p\equiv \pm 1 \bmod{8}$ and $O(5,p,0)$ if $p \equiv \pm 3 \bmod{8}$. 
These are examples of {\em locally toroidal\/} regular polytopes, meaning 
that their facets and vertex-figures are spherical or toroidal (but not all 
spherical). Such polytopes have not been fully classified (see 
\cite[Ch.12]{arp}). The three closely related groups $[3,4,3,3,3]$, 
$[3,3,4,3,3]$ and $[4,3,3,4,3]$ of rank $6$ similarly give rise to locally 
toroidal $6$-polytopes with groups $O_1(6,p,+1)$ if $p\equiv \pm 1 \bmod{8}$ 
and $O(6,p,+1)$ if $p \equiv \pm 3 \bmod{8}$. Each polytope is covered by 
the respective universal polytope 
\[\begin{array}{c}
\{ \,\{3,4,3,3\}_{(p,0,0,0)} \, , \, \{4,3,3,3\}\, \}, \\
\{ \,\{4,3,3,4\}_{(p,0,0,0)} \, , \, \{3,3,4,3\}_{(p,0,0,0)}\, \}, \\
\{ \,\{3,3,4,3\}_{(p,0,0,0)} \, , \, \{3,4,3,3\}_{(p,0,0,0)}\, \} ,
\end{array}\]
which has been conjectured to exist for all primes $p\geq 3$ and to be 
infinite for $p>3$ (see \cite[Ch.12]{arp}). 
For $p=3$ (and $p=2$, but 
this is outside our discussion), these three universal polytopes are known 
to be finite. The first two have group $\mathbb{Z}_3 \rtimes O(6,3,+1)$, and 
the last has group $(\mathbb{Z}_3 \oplus \mathbb{Z}_3) \rtimes O(6,3,+1)$ 
(see \cite{monsch3}).

More generally we can reduce the four groups $G$ modulo any integer 
$d\geq 2$ to obtain other kinds of locally toroidal regular polytopes 
of ranks $5$ or $6$, now with type vectors for facets and vertex-figures 
depending on $\Delta(G)$ and its subdiagrams for facets and vertex-figures,
as well as on the parity of $d$ (see \cite{monschmod}). In particular, 
this establishes \cite[Conjecture 12C2]{arp} concerning the existence of 
locally toroidal regular $6$-polytopes of type $\{3,4,3,3,3\}$, saying 
that the universal regular $6$-polytopes 
$\{\{3,4,3,3\}_{(d,0,0,0)},\{4,3,3,3\}\}$ and
$\{\{3,4,3,3\}_{(d,d,0,0)},\{4,3,3,3\}\}$ exist for all $d\geq 2$; this then 
settles 
the existence of the first kind of $6$-polytopes mentioned in the previous 
paragraph. 
The corresponding conjectures for the locally toroidal regular 
$6$-polytopes of types $\{3,3,4,3,3\}$ and $\{4,3,3,4,3\}$ are still 
open (see \cite[12D3,12E3]{arp}).

Another  interesting special class consists of 
the {\em $3$--infinity groups\/} 
$G=[p_1,\ldots,p_{n-1}]$, for which all periods $p_j \in \{3, \infty\}$. 
Typically then we have an alternating string of $3$'s and $\infty$'s, and 
the outcome depends on the nature of the string. For the prime $p=3$ we 
always have a C-group, namely $S_{n+1}$. 
If $G = [3^k, \infty^l]$, with 
$k + l = n-1$, then for all primes $p \geq 3$, the group $G^p$ is a 
string $C$-group. We have seen examples with $n=3$ or $4$. For 
$G = [3, \infty^l, 3]$, with $n= l+3, l\geq 1$, we have a string 
$C$-group, except possibly when $p = 7$ and $l \geq 4$, with 
$l \equiv 1 \bmod 3$. On the other hand, if 
$G = [3^k, \infty^l, 3^m]$, with $k>1$ or $m >1$, and $l\geq 1$, then 
$G^p$ is a string $C$-group for all but finitely many primes $p$. However, 
$G^p$ very often is not a C-group. For example, suppose that $G$ has a string
subgroup of the form $[\ldots , \infty, 3^k, \infty, \ldots]$, with 
$k \geq 1$, and that $p\geq 5$; then $G^p$ is not a string $C$-group, 
except possibly when $p=5$ or $7$ and $k \leq 1$ for all such string 
subgroups $[\ldots , \infty, 3^k, \infty, \ldots]$.

\section{Polytopes from hyperbolic Coxeter groups}
\label{polhyp} 

Following \cite{chir1} we begin with the symmetry group 
$\Gamma = [r,s,t]$ of a regular honeycomb $\{r,s,t\}$ in 
hyperbolic $3$-space $\mathbb{H}^3$ and consider a faithful 
representation of $\Gamma$ as a group of complex M\"obius 
transformations. Recall here that the absolute of $\mathbb{H}^3$ can 
be identified with the complex inversive plane (compare the Poincare 
halfspace model of $\mathbb{H}^3$), and that every group of hyperbolic 
isometries is isomorphic to a group of M\"obius transformations over 
$\mathbb{C}$. In particular, a plane reflection in $\mathbb{H}^3$ 
corresponds to an inversion 
in a circle (or line) determined as the 
``intersection" of the mirror plane in $\mathbb{H}^3$ with the absolute.
Under this correspondence the generating plane reflections 
$\rho_0,\ldots,\rho_3$ of the hyperbolic reflection group 
$\Gamma$ become inversions in circles, again denoted by 
$\rho_0,\ldots,\rho_3$, cutting one another at the same 
angles as the corresponding reflection planes in $\mathbb{H}^3$. 
The group $\Gamma$ is one of ten possible groups, namely
\[ [4,4,3],\, [4,4,4],\, [6,3,3],\, [6,3,4],\, [6,3,5],\,
 [6,3,6],\, [3,6,3],\, [3,5,3],\,[5,3,4],\, [5,3,5] ;\]
see \cite{weinc} for a complete list of the generating inversions 
for these groups.

Next recall that M\"obius transformations may conveniently be represented
(uniquely up to scalar multiplication) by $2\times 2$ matrices, namely
\[ \frac{az+b}{cz+d} \longleftrightarrow 
\left[ \begin{array}{cc}
a&b\\
c&d
\end{array} \right]
\qquad \mbox{ and }\qquad
\frac{a\overline{z}+b}{c\overline{z}+d} \longleftrightarrow \!\!\!\!\!\!
{\begin{array}{cc}
\\
\\
\end{array}}^{\#}\!\!\left[ \begin{array}{cc}
a&b\\
c&d
\end{array} \right], \]
with appropriate interpretations for multiplication of such matrices. We 
then obtain matrices $r_0,r_1,r_2,r_3$ and $s_1,s_2,s_3$, respectively, 
for the generators $\rho_0,\rho_1,\rho_2,\rho_3$ of $\Gamma$ and 
$\sigma_1,\sigma_2,\sigma_3$ of $\Gamma^{+}$. Let $G$ denote the group of 
complex $2\times 2$ matrices generated by $s_1,s_2,s_3$. When considered 
modulo scalars, $G$ (or rather the corresponding projective group $PG$) is 
isomorphic to $\Gamma^+$, the rotation group of $\{r,s,t\}$. 
It turns out that the matrices in $G$ all have entries in a certain subring 
$\mathbb{D}$ of $\mathbb{C}$ depending on $[r,s,t]$. Thus modular 
reduction applies. We
choose appropriate ideals $J$ in 
$\mathbb{D}$ and consider the matrices $s_1,s_2,s_3$ over the quotient ring 
$\mathbb{D}/J$, again modulo scalars (determined by a suitable subgroup of the 
center).  Then under certain conditions, the resulting projective group $G^J$ 
(say) is either the rotation subgroup for a directly regular polytope, or the full automorphism group of a chiral 
polytope. 

We illustrate the method for the hyperbolic group $[4,4,3]$. When 
$\Gamma = [4,4,3]$ is viewed as a group of complex M\"obius transformations, 
the generators  may be taken as
\[ \rho_{0}(z) = \overline{z},\quad \rho_{1}(z) = i\overline{z},\quad
\rho_{2}(z) = 1-\overline{z},\quad \rho_{3}(z) = 1/\overline{z}. \]
Then the rotation subgroup $\Gamma^+$ of $\Gamma$ is generated by 
\[ \sigma_{1}=\rho_{0}\rho_{1}=-iz, \quad 
\sigma_{2}=\rho_{1}\rho_{2}=-iz+i, \quad  
\sigma_{3}=\rho_{2}\rho_{3}=1- 1/z, \]
and consists only of proper M\"obius transformations (not involving complex 
conjugation). Here the matrices $s_1,s_2,s_3$ for $\sigma_1,\sigma_2,\sigma_3$ 
are given by    
\[ s_1 = \left[ \begin{array}{rr}
-i&0\\
0&1
\end{array} \right], \quad
s_2 = \left[ \begin{array}{rr}
-i&i\\
0&1
\end{array} \right], \quad
s_3 = \left[ \begin{array}{rr}
1&-1\\
1&0
\end{array} \right], \]
and $G$ is the group of all invertible $2\times 2$ matrices over the Gaussian 
integers $\mathbb{D}=\mathbb{Z}[i]$ (with determinants $\pm 1, \pm i$). In 
particular,
$[4,4,3]^{+} \cong PSL_2(\mathbb{Z}[i]) \rtimes C_2$,
where the first factor is also known as the Gaussian modular group or Picard 
group (see also \cite{joh}).
Since we wish to obtain toroidal facets of type $\{4,4\}^{+}_{(b,c)}$, we must 
consider the imposition of the extra relation 
$(\sigma_{1}^{-1}\sigma_2)^b (\sigma_1\sigma_2^{-1})^c = 1$
on the generators of $\Gamma^+$ (see \cite{genrel}). Note here that the 
corresponding matrix product in $G$ is
\[ (s_{1}^{-1} s_2)^b (s_1 s_2^{-1})^c 
\,=\, \left[ \begin{array}{rr}
1&-1\\
0&1
\end{array} \right]^b \,
\left[ \begin{array}{rr}
1&-i\\
0&1
\end{array} \right]^c
\,=\, \left[ \begin{array}{rc}
1&-(b+ci)\\
0&1
\end{array} \right]. \]
Thus the ideals $J$ of $\mathbb{Z}[i]$ must be chosen in such a way that 
$b+ci=0$ in $\mathbb{Z}[i]/J$.

One natural choice of ideal is $J =m\mathbb{Z}[i]$, where $m\geq 3$ is an 
integer in $\mathbb{Z}$. This choice of ideal typically produces directly regular polytopes of type 
$\{\{4,4\}_{(m,0)},\{4,3\}\}$ whose rotation subgroup is 
$PSL_2(\mathbb{Z}_{m}[i])$ or a closely related group, with the exact
structure determined by the prime factorization of $m$ in $\mathbb{Z}$ 
(see \cite[p.238]{chir1}). For example, if $m \equiv 3 \bmod 4$ is a 
prime, we obtain $PSL_2(\mathbb{Z}_{m}[i]) \cong PSL_2(m^2)$.

A more interesting choice of ideal arises from a solution of the equation 
$x^{2}=-1 \bmod m$, with $m$ as above. Let 
$m= 2^{e}p_{1}^{e_1}\ldots p_{k}^{e_k}$ be the prime factorization in 
$\mathbb{Z}$. Then the equation is solvable if and only if $e=0,1$and $p_{j}\equiv 1\bmod 4$ for each $j$. If $\hat{i}\in \mathbb{Z}_m$ is
such that $\hat{i}^{2}=-1 \bmod m$, then there exists a unique pair of 
positive integers $b,c$ such that $m=b^2+c^2$, $(b,c)=1$ and 
$b = -\hat{i}c \bmod m$. We now take $J = (b+ci)\mathbb{Z}[i]$, 
which is  the kernel of the ring epimorphism 
$\mathbb{Z}[i] \rightarrow \mathbb{Z}_m$ that maps the complex number 
$x+yi$ to the element $x_m+y_m\hat{i}$ in $\mathbb{Z}_m$, where 
$x_m\equiv x$ and $y_m\equiv y \bmod m$. This choice of ideal  ``destroys" the reflections in the overlying reflection group and 
typically yields chiral polytopes of type $\{\{4,4\}_{(b,c)},\{4,3\}\}$ 
whose rotation subgroup is $PSL_2(\mathbb{Z}_{m})$ or a closely related 
group, with the exact structure again determined by the prime factorization 
of $m$. For example, if $m \equiv 1 \bmod 8$ is a prime, then the group is
$PSL_2(m)$. Note that the structure of the polytope also depends on the 
solution $\hat{i}$.  
For example, when $m=65$, the solution $\hat{i}=8$ 
gives facets $\{4,4\}_{(1,8)}$, while $\hat{i}=18$ leads to facets 
$\{4,4\}_{(4,7)}$.

With similar techniques we can construct a host of regular or chiral
 $4$-polytopes
of types $\{4,4,4\}$, $\{4,4,4\}$, $\{6,3,3\}$, $\{6,3,4\}$, 
$\{6,3,5\}$, $\{6,3,6\}$ or $\{3,6,3\}$ (see \cite{nosc,chir1}). 
However, the underlying ring $\mathbb{D}$ will depend on the Schl\"afli 
symbol. The Gaussian integers $\mathbb{Z}[i]$ and Eisenstein integers 
$\mathbb{Z}[\omega]$ (with $\omega$ a cube root of unity) suffice in 
all cases except $\{6,3,4\}$ and $\{6,3,5\}$ (this is based on subgroup 
relationships between the rotation subgroups for the various types). For 
$\{6,3,4\}$ and $\{6,3,5\}$, respectively, we can work over the ring 
$\mathbb{Z}[\omega,\sqrt{2}]$ and $\mathbb{Z}[\omega,\tau]$ (with $\tau$ the 
golden ratio; see the next section). 

Chiral polytopes also exist in ranks larger than $4$, but explicit 
constructions of finite examples (or rank $5$) were only discovered 
quite recently in \cite{co4}. Modular reduction is an effective method to 
produce examples in rank $4$. In addition to the symmetry groups of 
$3$-dimensional hyperbolic honeycombs, several other discrete hyperbolic
groups admit representations as groups of linear fractional transformations 
over other rings of complex or quaternionic integers (described in detail in 
the forthcoming book by Johnson~\cite{joh}).  Although the arithmetic involved 
is likely to be considerably more complicated than in the rank $4$ case, 
there is a good chance that the reduction method will carry over to produce 
examples of finite chiral polytopes of rank larger than $4$.

\section{The groups $[3,5,3]$ and $[5,3,5]$}
\label{noncrys}
\vskip.05in\noindent

In Section~\ref{polcry} we observed that any crystallographic
string Coxeter group $G$ has
rotational periods $p_j \in \{ 2,3,4,6,\infty\}$ and can be represented 
faithfully 
as a matrix group over the domain $\mathbb{Z}$. Here we widen the 
discussion a little by allowing  the period $p_j = 5$.
Keeping (\ref{stanfm}) in mind, we note that 
$2\cos\frac{\pi}{5} = \tau$,
where the golden ratio  $\displaystyle{\tau =  (1+\sqrt{5})/2}$
is  the positive root of
$\tau^2 = \tau +1$. We therefore move to the larger coefficient domain
$$ \mathbb{D}:= \mathbb{Z}[\tau] = \{ a+b\tau\, : \, a,b \in \mathbb{Z} \}\;, $$ 
and soon find that we need only  add the subdiagram
$$\stackrel{s}{\bullet}\!\!\!-\!\!\!-\!\!\!-\!\!
\!\!\stackrel{\tau^2 s}{\bullet}$$ 
to those already listed in  
Table~\ref{binartabII} in order for the Cartan integers $m_{ij}$ of (\ref{reflnII}) to be in $\mathbb{D}$ for all $i,j$. This subdiagram, say on nodes $i,j$, 
does indeed define the non-crystallographic dihedral group
$\langle r_i, r_j \rangle$  with order $10$ and period $p_{ij} = 5$.
(In other notation, this is the group $H_2 \simeq I_2(5)$.)
Naturally, we must now allow rescaling of nodes by any `integer'
$s \in \mathbb{D}$ or its inverse.  Furthermore, referring back to 
(\ref{reflnII}), we find that
$m_{ij} = \tau^2 \in \mathbb{D}$, so that $G$ is represented
as a matrix group over $\mathbb{D}$ through its action on 
the $\mathbb{D}$-module 
 ${\displaystyle \oplus_{j} \mathbb{D} b_j}$.

Let us now summarize the key  arithmetic properties of the domain $\mathbb{D}$.
(We refer to \cite{dodd} for a detailed account of this,
and to \cite{chen} for a deeper discussion of `$\mathbb{D}$-lattices'
and the related finite Coxeter groups $H_k,\, k = 2,3,4$.)
First of all, we recall that $\mathbb{D}$
is the ring of algebraic integers in 
the field $\mathbb{Q}(\sqrt{5})$. 
The non-trivial field automorphism
mapping $\sqrt{5} \mapsto -\sqrt{5}$ induces a ring  automorphism
$' : \mathbb{D} \rightarrow \mathbb{D}$, 
which in this section we shall call \textit{conjugation}. Thus
$$(a+b\tau)' = (a+b)-b\tau\;.$$
In particular,
${\tau}' = 1-\tau =  -\tau^{-1}$. Recall that $z = a+b\tau$
has  \textit{norm} $N(z) := z  z' =  a^2 + ab -b^2$.
We note   that $\mathbb{D}$ is a Euclidean domain, through a  division
algorithm based on   $|N(z)|$.

The set of \textit{units} in $\mathbb{D}$ is
$ \{ \pm \tau^n\, : \, n \in \mathbb{Z}\} \; = \; \{ u \in \mathbb{D}\, : \,
N(u) = \pm 1\}.$
Recall that integers $z, w \in \mathbb{D}$ are \textit{associates} if 
$z = u w$ for some unit $u$. Up to associates,
the primes $\pi \in \mathbb{D}$ can be classified as follows:

\begin{itemize}
\item the prime $ \pi = \sqrt{5} = 2\tau -1$,
which is self-conjugate (up to associates: $\pi' = -\pi$);
\item  rational primes $ \pi = p \equiv \pm 2 \bmod{5}$, also self-conjugate; 
\item  primes $\pi = a+b\tau$, for which   
$|N(\pi)|$ equals a rational prime
 $ q \equiv \pm 1 \bmod{5}$. In this case,  the conjugate prime  
 $\pi' = (a+b)-b\tau$ is  not an associate of $\pi$.  
\end{itemize}

Let us now turn to the group $G = [3,5,3]$, here acting as an
orthogonal group on real $4$-space $V$ (in contrast to the conformal 
representation on $\mathbb{H}^3$ mentioned in the previous section).
Since $\tau^2$ is a unit, there is essentially only one choice of 
diagram, namely 
 
$$ \Delta(G) = \;\;\stackrel{1}{\bullet}\frac{}{\;\;\;\;\;\;\;}
	\stackrel{1}{\bullet}\frac{}{\;\;\;\;\;\;\;}
	\stackrel{\tau^2}{\bullet}\frac{}{\;\;\;\;\;\;\;}
	\stackrel{\tau^2}{\bullet} \;.
$$
The discriminant is 
$$\mbox{disc}(V) = -\frac{1}{16} (2+5\tau) \sim -(2+5\tau)\;,$$  
where the prime $\delta := -(2+5\tau)$ has norm $-11$.

Now consider any prime $\pi \in \mathbb{D}$. Our goals are to show that 
$G^\pi = \langle r_0, r_1, r_2, r_3 \rangle^\pi$ is a  
string $C$-group and to determine its structure, then
say a little about the corresponding polytope 
$\mathcal{P}^\pi := \mathcal{P}(G^\pi)$. 

In fact, we can almost immediately apply a suitable 
generalization of Theorem 4.2  in \cite{monsch1}. First note that the subgroup 
$G_3^\pi = \langle r_0, r_1, r_2 \rangle^\pi$ is obviously some quotient
of the spherical group $[3,5] \simeq H_3 $. Now it is easy to check that
after reduction modulo any prime $\pi$, even for associates of $2$, the
reflections $r_j$ still have period $2$. Next, we consider the isometry
\begin{equation}\label{facpet}
z:= (r_0 r_1 r_2)^5 = \left[
\begin{array}{cccc}
  -1& 0&0 & \tau^4\\
  0& -1&0 & 2\tau^4\\
 0 & 0& -1& 3\tau^2\\
  0& 0& 0& 1\\
\end{array}
\right]\;\;.
\end{equation}
in $G$. Since $\tau^4$ is a unit, this means that $r_0 r_1 r_2$ still has period 
$10$ in $G^\pi$. Thus $\langle r_0, r_1, r_2 \rangle^\pi \simeq [3,5]$
and dually $\langle r_1, r_2, r_3 \rangle^\pi \simeq [5,3]$.
Consulting the proof of \cite[Th. 4.2]{monsch1}, we see that we
 need only show that the 
orbit of $\mu_0  := [1,0,0,0]$  under the right action of the matrix
group $\langle r_0, r_1, r_2 \rangle$ has the same size modulo $\pi$
as in characteristic $0$, namely $12$. This is routinely verified, so we
have proved most of

\begin{PROP} Let $G=[3,5,3]$. For any prime $\pi \in \mathbb{D}$, the group
$G^\pi = \langle r_0, r_1, r_2, r_3 \rangle^\pi$ is a finite string 
$C$-group. The corresponding finite regular polytope $\mathcal{P}^\pi$
is self-dual and has
icosahedral facets $\{3,5\}$ and dodecahedral vertex figures $\{5,3\}$.
\end{PROP}

\noindent
\textbf{Proof}.  To verify self-duality we define $g \in GL(V)$
by $g: [b_0, b_1, b_2, b_3] \mapsto [\tau^{-1} b_3, \tau^{-1} b_2, \tau b_1, \tau b_0]$. 
Then  $g^2 = 1$, 
$g r_0 g = r_3$ and $ g r_1 g =  r_2$.  (See \cite[2E12]{arp}.)
\hfill$\square$

\medskip

A more detailed description of $G^\pi$ must depend on 
the nature of the prime $\pi$. (Of course, our results
are typically  unaffected by replacing $\pi$ by any associate
$\pm \tau^m \pi$.)
In all cases the underlying finite field
$\mathbb{K} := \mathbb{D}/(\pi)$ has order $|N(\pi)|$, so that
$G^\pi$ acts as an orthogonal group on the $4$-dimensional
vector space $V$ over $\mathbb{K}$ preserving the modular image of 
the bilinear form for $G$.

\noindent\textbf{Case 1}: $\pi = 2$.

Here an easy  calculation using  GAP
confirms that $G^\pi$ is the orthogonal group
$O(4, 2^2,-1)$ with Witt index $1$ over $\mathbb{K} = GF(2^2)$.
Since $|G^2| = 8160$, the polytope $\mathcal{P}^2$ has $68$ vertices
and $68$ icosahedral facets.

\medskip

Henceforth we suppose that $\pi$ is not an associate of $2$.
To work with such primes, we need a generalization of the rational 
Legendre symbol $\Leg{p}{q}$. Thus for any $\alpha \in \mathbb{D}$ 
and  prime $\pi$ we set

$$
\LegD{\alpha}{\pi} :=\left\{ 
\begin{array}{cl}
+1, & \mbox{ if $\alpha$ is a quadratic residue (mod $\pi$);}\\
-1, & \mbox{ otherwise.}
\end{array}
\right.
$$
(Compare \cite[Ch. VIII]{dodd}.)
We are mainly interested in computing
$$\epsilon: = \LegD{\delta}{\pi}$$
where $\delta = -(2 + 5\tau)$ is the discriminant.
Since every label in $\Delta(G)$ is square, we conclude that $G^\pi$ is a
subgroup of $O_1(4, |N(\pi)|,\epsilon)$,
so long as $\delta$ and $\pi$ are relatively prime. Indeed, $G^\pi$
will almost  always equal such an orthogonal group.

\medskip
\noindent\textbf{Case 2}: $\pi = \sqrt{5} = 2\tau-1$. 

Here $|\pi \pi'| = 5 \equiv  0\bmod{\pi}$, so that the discriminant
$\delta = -(2+5\tau) \equiv 3 \bmod{\pi}$, which
is non-square in $\mathbb{K} = GF(5)$. Thus $\epsilon = -1$
and $G^\pi = O_1(4,5,-1)$ has order $15600$. In fact, the polytope 
$\mathcal{P}^{\sqrt{5}}$ is isomorphic to that obtained in 
\cite[p. 347]{monsch2}
through reduction $\bmod{5}$ of the crystallographic group $[3, \infty, 3]$.

\medskip
\noindent\textbf{Case 3}: $\pi$ is an associate of an odd rational
prime $p \equiv \pm 2 \bmod{5}$.

Since $\mathbb{K} = GF(p^2)$, $G^\pi = G^p$ is, for suitable $\epsilon$, 
a subgroup
of $O_1(4,p^2, \epsilon)$, whose order we recall
is $p^4 (p^4 - \epsilon) (p^4 -1)$.
Consulting  \cite[Th. 3.1]{monsch1},
we see that
$G^p = O_1(4,p^2, \epsilon)$ so long as we can rule out two remote 
alternatives.

First of all, it is conceivable that $G^p \simeq H_4 = [3,3,5]$. But here
it is easy to check directly that $H_4$ cannot be generated by reflections 
$r_j$  satisfying the Coxeter-type relations inherited from $[3,5,3]$,
let alone the independent relations induced by reduction modulo ${p}$.

Secondly, we must show that $G^p$ is not isomorphic to some orthogonal group
$O_1(4,p,\eta)$, $\eta = \pm 1$, over the subfield $GF(p)$. If this were so, 
then Theorem 3.1 in \cite{monsch1} would actually imply that
$G^p$ is similar to $O_1(4,p,\eta)$  under extension of scalars. More precisely,
if $\mathbb{L}$ is an algebraic closure of $\mathbb{K} = GF(p^2)$, 
then there would
exist some $g \in GL(V_{\mathbb{L}})$ with $g G^p g^{-1} = O_1(4,p,\eta)$.
Using (\ref{reflnII}), we   compute with respect to the new
basis $\{c_i\} = \{g(b_i)\}$
for $V_{\mathbb{L}}$. Thus the reflection
$\tilde{r_i} := g r_i g^{-1}$ satisfies
$$ \tilde{r_i}(c_j) = g(b_j + m_{i j} b_i) = c_j + m_{i j} c_i\;.
$$
We conclude that  the field of definition for $G^p$ must always 
contain the subfield 
generated by the Cartan integers $m_{i j}$. In our case, 
$m_{12} = \tau^2 \not\in GF(p)$, so that $g G^p g^{-1}$ cannot possibly
be a group $O_1(4,p,\eta)$.

Having shown that $G^p = O_1(4,p^2,\epsilon)$, we next
 determine $\epsilon$. From \cite[Th. 8.5(a)]{dodd} we have
$$\LegD{\delta}{\pi} = \Leg{N(\delta)}{p} = \Leg{-11}{p}
=\Leg{-1}{p} \Leg{11}{p} = \Leg{p}{11}\; ,
$$
by (rational) quadratic reciprocity. Since the non-zero 
squares $(\mbox{mod }{11})$
are $1,3,4,5,9$, we conclude that
$$\epsilon :=\left\{ 
\begin{array}{cl}
+1, & \mbox{ if }  p \equiv 3, 12, 23, 27, 37, 
38, 42, 47, 48, 53 \bmod{55};\\
-1, & \mbox{ if } p \equiv 2, 7, 8, 13, 17, 18, 28, 32, 43, 52 \bmod{55}.
\end{array}
\right.
$$

\medskip
\noindent\textbf{Case 4}: $ \pi = a+ b\tau$, where 
$N(\pi) = a^2 + ab -b^2 = q$, where the rational prime  
$q \equiv \pm 1\bmod{5}$; however, $\pi$ is 
not an associate of $\delta = -(2+5\tau)$.

We now have $\mathbb{K} = GF(q)$. An even easier appeal to 
\cite[Th. 3.1]{monsch1} gives $G^\pi = O_1(4,q,\epsilon)$. We need only 
determine $\epsilon = \LegD{\delta}{\pi}$.
Since $a+b\tau \equiv 0 \bmod{\pi}$, we may suppose 
$\tau = -a/b \in \mathbb{K}$. Thus
$$ \delta = -(2+5\tau) \sim -b^2(2+5\tau) \equiv 5ab -2b^2 \bmod{\pi}\;.
$$
By \cite[Th. 8.5(a)]{dodd}, we obtain
$$ \epsilon = \LegD{\delta}{\pi} = \Leg{(5ab -2b^2)}{q} = 
\Leg{b}{q} \Leg{(5a -2b)}{q}\;.$$
Using the rational Legendre symbol, we can thus compute $\epsilon$
for any prime $\pi = a+b\tau$.

It is possible to say when $\pi$ and its conjugate $\pi'$ 
give opposite $\epsilon$'s, so that the corresponding orthogonal
spaces have, in some order,
Witt indices $1$ and $2$. This happens if and only if
$q$ is a square mod $11$, 
since
$$ \LegD{\delta}{\pi} \, \LegD{\delta}{\pi'} = 
\Leg{(b(5a-2b)(-b)(5a+7b))}{q} = \Leg{-11}{q} =\Leg{q}{11}\;.$$

One notable instance here is $\pi = \delta' = -7 + 5\tau$, which is 
relatively prime to  the discriminant $\delta$. We have
$G^{\delta'} = O_1(4,11,-1)$, of order 1771440.

\medskip
\noindent\textbf{Case 5}: $\pi = \delta = -(2+5\tau)$.

This is the only case in which the orthogonal space $V$ is singular. 
Now $\mathbb{K}=GF(11)$ and $\tau = -2/5 = 4$. 
We find that   $\mbox{rad}(V)$ is spanned by 
$c = 7 b_0 + 3 b_1 +2b_2 + b_3$, and that
$V = \mbox{rad}(V) \perp V_3$, where $V_3$ is the non-singular 
subspace spanned by $b_0, b_1, b_2$. It is then not hard to see that
$$O(V) \simeq \check{V_3} \rtimes (\mathbb{K}^* \times O(V_3))\;,$$
where $\mathbb{K}^* \simeq GL(\mbox{rad}(V))$ and 
$\check{V_3}$ is dual to $V_3$. We observe that the abelian
group $\check{V_3} \simeq \mathbb{K}^3$  consists of
all transvections 
$$r(x) = x+ \varphi(x) c\;,$$
where $\varphi \in \check{V_3}$ (with $\check{V_3}$ viewed as a 
subspace of $\check{V}$ fixing $c$). Now since every $r_j$ fixes $c$, 
$G^\delta$ must be  a subgroup of  the pointwise stabilizer of $\mbox{rad}(V)$.
In fact, another calculation with GAP  confirms that

$$G^\delta = \widehat{O}_1(V) \simeq  
\check{V_3} \rtimes  O_1(V_3) \;,$$
which has order $11^3 \cdot 11 \cdot (11^2 -1) = 1756920$. Now consider the 
isometry $z\in G$ defined in (\ref{facpet}). It is easy to check that
$z(c)\equiv c \bmod{\delta}$, so that 
$z = 1_{\rm{rad}(V)}  \perp  -1_{V_3} \in \widehat{O}_1(V)$ acts as
the  \textit{central inversion} in the group $O_1(V_3)$ 
for the icosahedral facet. Thus $G^\delta$ has a normal
subgroup $A$ isomorphic to 
$\check{V_3} \rtimes \langle z\rangle$
and so of order $2 \cdot 11^3$. Using 
$O_{1}(3,11,0) \simeq PSL_2(11) \rtimes C_2$ (see \cite[Th. 5.20]{art}), 
we conclude that
$$\overline{G} := G^{\delta}/A  \simeq PSL_2(11)\; ,$$
of order $660$.
Remarkably, $\overline{G}$ is also a string $C$-group. The resulting
polytope is the $11$-cell independently
discovered by Coxeter in \cite{11cell} and Gr\"{u}nbaum in \cite{Gru}.
Indeed, both $r_0 r_1 r_2$ and $r_1 r_2 r_3$ have period $5$ in the quotient   
(see (\ref{facpet})), and 
$$\mathcal{P}(\overline{G}) = \{\, \{3,5\}_5\, , \, \{5,3\}_5\,\}$$
is the universal $4$-polytope with hemi-icosahedral facets 
and hemi-dodecahedral vertex-figures.

This finishes our investigation of the group $[3,5,3]$. 
Evidently a somewhat similar  analysis is possible for
the group $H = [5,3,5]$ with diagram
$$ \Delta(H) = \;\;\stackrel{1}{\bullet}\frac{}{\;\;\;\;\;\;\;}
	\stackrel{\tau^2}{\bullet}\frac{}{\;\;\;\;\;\;\;}
	\stackrel{\tau^2}{\bullet}\frac{}{\;\;\;\;\;\;\;}
	\stackrel{1}{\bullet}  
$$
and corresponding  discriminant 
$\frac{-1}{16}(3+7\tau) \sim -(3+7\tau) =: \lambda$. Since  
$N(\lambda) = -19$, we see that $\lambda$ is also prime. 
We note only that the group $H^\lambda$ for the singular space $V$
again has an  interesting quotient. In fact,
$$\overline{H} \simeq PSL_2(19)$$
is the automorphism group for the universal regular polytope
$$\mathcal{P}(\overline{H}) = \{\, \{5,3\}_5\, , \, \{3,5\}_5\,\}\;,$$ 
with hemi-dodecahedral facets and hemi-icosahedral vertex-figures. This is the $57$-cell described by Coxeter in \cite{57cell}.

With the exception of the $11$-cell and $57$-cell, the polytopes described here can be viewed as regular tessellations on hyperbolic $3$-manifolds (see \cite[6J]{arp}). Moreover, the two exceptions are the only regular polytopes or rank $4$ (or higher) with automorphism group isomorphic to $PSL_2(r)$ for some prime power $r$ (see \cite{lee}). 
For related work see also \cite{hal1,hal2,jones}.
\\[-.4in]

\medskip
\noindent
{\em B. Monson, University of New Brunswick, Fredericton, New Brunswick, Canada E3B 5A3, bmonson@unb.ca\/}\\[.1in]
\noindent
{\em Egon Schulte, Northeastern University, 
Boston, Massachussetts,  USA, 02115, schulte@neu.edu\/}

\end{document}